\documentclass[10pt]{amsart}
\usepackage[T1]{fontenc}
\usepackage[latin5]{inputenc}
\usepackage{float}
\usepackage{amssymb}
\usepackage{amsfonts}
\usepackage[centertags]{amsmath}
\usepackage{amsthm}
\usepackage{newlfont}
\usepackage{array}
 \makeatletter
\theoremstyle{plain}
\newtheorem{thm}{Theorem}[section]
\numberwithin{equation}{section}
\numberwithin{figure}{section}  
\theoremstyle{plain}

\theoremstyle{plain}

\theoremstyle{plain}
\newtheorem{cor}[thm]{Corollary} 
\theoremstyle{plain}

\theoremstyle{plain}
\newtheorem{lem}[thm]{Lemma} 
\theoremstyle{plain}
\newtheorem{rem}[thm]{Remark}

\begin{document}
\title{A note on the some geometric properties of the sequence spaces defined by Taylor method}
\author{Murat Kiri\c{s}ci}
\address[Murat Kiri\c{s}ci]{Department of Mathematical Education, Hasan Ali Y\"{u}cel Education Faculty, Istanbul University, Vefa, 34470, Fatih, Istanbul, Turkey}
\email{mkirisci@hotmail.com, murat.kirisci@istanbul.edu.tr}
\vspace{0.5cm}
\vspace{0.5cm}

\begin{abstract}
In this paper, it was obtained the new matrix domain with the well known classical sequence spaces and an infinite matrix.
The Taylor method which known then as the circle method of order $r$ (0 < r < 1), as an infinite matrix for the matrix domain is used.
Newly constructed space is isomorphic copy of the spaces of all absolutely $p-$summable sequences. It is well known that Hilbert space 
have the nicest geometric properties. Then, it is proved that the new space is a Hilbert space for $p=2$. Further, it was computed 
dual spaces and characterized some matrix classes of the new Taylor space in the table form. Section 3 is devoted some geometric properties of Taylor space.
\end{abstract}

\maketitle

\section{Introduction}

The $e_{p}^{r}$ and $e_{\infty}^{r}$ sequence spaces using Euler mean were defined by Altay et al. \cite{abm}, as follows:
\begin{eqnarray*}
e_{p}^r&:=&\left\{x=(x_{k})\in \omega: \sum_{n} \left|\sum_{k=0}^n \binom{n}{k}(1-r)^{n-k}r^{k}x_{k}\right|^{p}<\infty\right\}, \quad (1\leq p < \infty),\\
e_{\infty}^r&:=&\left\{x=(x_{k})\in \omega: \sup_{n\in\mathbb{N}}\left|\sum_{k=0}^n \binom{n}{k}(1-r)^{n-k}r^{k}x_{k}\right|<\infty\right\},
\end{eqnarray*}
where $E^{r}=(e_{nk}^r)$ denotes the Euler means of order $r$ defined by
\begin{eqnarray*}
e_{nk}^r = \left\{ \begin{array}{ccl}
\binom{n}{k}(1-r)^{n-k}r^{k}&, & (0\leq k\leq n),\\
0&, & (k > n)
\end{array} \right.
\end{eqnarray*}
for all $k,n\in\mathbb{N}$. It is known that the method $E^{r}$ is regular for $0<r<1$ and $E^{r}$ is invertible such that $(E^{r})^{-1}=E^{1/r}$ with $r\neq 0$. We assume unless stated otherwise that $0<r<1$.\\

Following Altay et al. \cite{abm} and  Mursaleen et al.\cite{mba}, Kirisci \cite{kiris2} defined the sequence spaces $t_{c}^{r}$ and $t_{0}^{r}$
consisting of all sequences $x=(x_{k})$ such that their Taylor transform $T(r)$, as below:

\begin{eqnarray*}
t_{0}^r&:=&\left\{x=(x_{k})\in \omega: \lim_{n \to \infty}\sum_{k=n}^{\infty} \binom{k}{n}(1-r)^{n+1}r^{k-n}x_{k}=0\right\},\\
t_{c}^r&:=&\left\{x=(x_{k})\in \omega: \lim_{n \to \infty}\sum_{k=n}^{\infty} \binom{k}{n}(1-r)^{n+1}r^{k-n}x_{k} \quad \textrm{exists}\right\}.\\
\end{eqnarray*}

Let $r\in \mathbb{C}/\{0\}$. Then, the Taylor transform $T(r)$ was defined by the matrix $T(r)=(t_{nk}^r)$, where
\begin{eqnarray}\label{mtrxTaylor}
t_{nk}^r = \left\{ \begin{array}{ccl}
\binom{k}{n}(1-r)^{n+1}r^{k-n}x_{k}&, & (k\geq n),\\
0&, & (0\leq k < n)
\end{array} \right.
\end{eqnarray}
for all $k,n\in\mathbb{N}$. In case $r=0$, it is immediate that $T(0)=(t_{nk}^{r})=(\delta_{nk})=I$.
In working with the Taylor method it is frequently used the fact that
\begin{eqnarray*}
\frac{1}{(1-z)^{n+1}}=\sum_{k=n}^{\infty}\binom{k}{n}z^{k-n}  ~\textrm{ if }~ |z|<1
\end{eqnarray*}
and additionally $\sum_{k=n}^{\infty}\binom{k}{n}z^{k-n}$ converges only if $|z|<1$.
The Taylor transform $T(r)$ is regular if and only if $0\leq r<1$, i.e., $r$ is real and $0\leq r<1$.
The product of the $T(r)$ matrix with the $T(s)$ matrix is the transpose of the matrix $(1-r)(1-s)E^{(1-r)(1-s)}$, where $E$ denotes Euler mean.
It is known that $T(r)$ is invertible and $(T(r))^{-1}=T[-r/(1-r)]$ provided $r\neq 1$.\\

The $\alpha-$(Köthe Toeplitz), $\beta-$(Generalized Köthe Toeplitz) and $\gamma-$(Garling)dual spaces
 of the sequence spaces $t_{c}^{r}$ and $t_{0}^{r}$ were computed and
necessary and sufficient conditions for the charaterization
of the matrix classess of $(t_{c}^{r}:Y)$, $(t_{0}^{r}:Y)$ and $(Y:t_{c}^{r})$, $(Y:t_{0}^{r})$ were given in \cite{kiris2}, where
$Y$ is one of the classical sequence space. Also, in \cite{kiris2}, Steinhaus type theorems were proved.\\

The aim of this paper is  examined some geometric properties of new spaces $t_{p}^{r}$ and $t_{\infty}^{r}$ defined by Taylor transform such as Banach-Saks,
weak Banach-Saks, weak fixed point, modulus of convexity.


\section{New Taylor Sequence Spaces of non-absolute type}

In this section, we give the definitions of the sequence spaces $t_{p}^{r}$ and $t_{\infty}^{r}$ and study some properties.\\

Following Altay et al.\cite{abm}, Mursaleen et al. \cite{mba} and Kirişci \cite{kiris2},
 we define the sequence spaces $t_{p}^r$ and $t_{\infty}^r$, as the set of all sequences such that $T(r)-$transforms of them are in the
spaces $\ell_{p}$ and $\ell_{\infty}$, that is
\begin{eqnarray*}
t_{p}^r&=&\left\{x=(x_{k})\in \omega: \sum_{n}\left|\sum_{k=n}^{\infty} \binom{k}{n}(1-r)^{n+1}r^{k-n}x_{k}\right|^{p}<\infty\right\}, \quad (1\leq p < \infty)\\
t_{\infty}^r&=&\left\{x=(x_{k})\in \omega: \sup_{n\in \mathbb{N}}\left|\sum_{k=n}^{\infty} \binom{k}{n}(1-r)^{n+1}r^{k-n}x_{k}\right|<\infty \right\}.\\
\end{eqnarray*}

If $\lambda$ is an arbitrary normed sequence space, then the matrix
domain $\lambda_{T(r)}$ is called \emph{Taylor sequence space}. This section of the study is a natural continuation of the results in \cite{kiris2}.
So, some topological properties of the spaces $t_{p}^r$ and $t_{\infty}^r$ are similar to the spaces $t_{c}^r$ and $t_{0}^r$.
 That is, in point of structural properties, we can obtain similar topological results
between new spaces and the spaces $t_{c}^r$ and $t_{0}^r$. Therefore, we won't give the proof of topological properties of new spaces. However, for $p=2$,
we show that the space $t_{p}^r$ is a Hilbert space. Also, we will prove the theorems of dual spaces which are important. \\

Define the sequence $y=\{y_{k}(r)\}$ by the $T(r)-$transform of a sequence $x=(x_{k})$, i.e.,
\begin{eqnarray}\label{transseq}
y_{k}(r)=\sum_{k=n}^{\infty}\binom{k}{n}(1-r)^{n+1}r^{k-n}x_{k} ~\textrm{ for all }~k\in\mathbb{N}.
\end{eqnarray}

The absolute property does not hold on the spaces $t_{p}^r$ and $t_{\infty}^r$, that is $\|x\|_{t_{p}^{r}} \neq \| |x| \|_{t_{p}^{r}}$
 and $\|x\|_{t_{\infty}^{r}} \neq \| |x| \|_{t_{\infty}^{r}}$ for at least one sequence in the spaces $t_{p}^r$ and $t_{\infty}^r$, where $|x|=(|x|)$.
This means that these spaces are sequence spaces of nonabsolute type.\\

The spaces $t_{p}^r$ and $t_{\infty}^r$ are linear spaces with coordinatewise addition and scalar multiplication that are $BK-$spaces
with the norm $\|x\|_{t_{p}^{r}}=\|T(r)x\|_{\ell_{p}}$, where $1\leq p \leq \infty$.\\

The space $t_{p}^r$ is an isomorphic copy of the classical $\ell_{p}$ space, for $1\leq p \leq \infty$ i.e., $t_{p}^r \cong \ell_{p}$ and $t_{\infty}^r \cong \ell_{\infty}$.

\begin{thm}
Except the case $p=2$, the space $t_{p}^r$ is not an inner product space, hence not a Hilbert space for $1\leq p < \infty$.
\end{thm}

\begin{proof}
For $p=2$, we will show that the space $t_{2}^r$ is a Hilbert space. It is known that the space $t_{p}^r$  for $1\leq p \leq \infty$, then the space $t_{2}^r$ is a $BK$-space, for $p=2$. Also its norm can be obtained from an inner product, i.e., $\|x\|_{t_{2}^r}=\langle T(r)x, T(r)x\rangle^{1/2}$ holds.
Then the space $t_{2}^r$ is a Hilbert space.\\

If we choose the sequences $e^{(0)}=(1,0,0,\cdots)$ and $e^{(1)}=(0,1,0,\cdots)$, we get that the norm of the space
 $t_{p}^r$ does not satisfy the parallelogram equality, which means that the norm cannot be obtained from inner product.
Hence, the space $t_{p}^r$ with $p\neq 2$ is a Banach space that is not a Hilbert space.
\end{proof}

Since the space $t_{p}^r$ is an isomorphic copy of the classical $\ell_{p}$ space, for $1\leq p \leq \infty$, then, we can consider the transformation $S$ from
$t_{p}^r$ to $\ell_{p}$ by $y=Sx=T(r)x$. Clearly, the transformation $S$ is a linear, surjective and norm preserving(From Theorem 2.2 of \cite{kiris2}). Because of the isomorphism $S$ is onto the inverse image of the basis of the space $\ell_{p}$ is the basis of the new spaces $t_{p}^{r}$. Therefore, choose a sequence $b^{(k)}(r)=\big\{ b_{n}^{(k)}(r) \big\}_{n\in \mathbb{N}}$ of elements of the space $t_{p}^r$ for every fixed $k\in \mathbb{N}$ by

\begin{eqnarray*}
b_{n}^{(k)}(r)= \left\{ \begin{array}{ccl}
\binom{k}{n}(1-r)^{-(k+1)}(-r)^{k-n}&, & \quad k\geq n\\
0&, & 0\leq k <n
\end{array} \right.
\end{eqnarray*}
and let $\lambda_{k}(r)=(T(r)x)_{k}$ for all $k\in \mathbb{N}$. Then,
 the sequence $\bigl\lbrace b^{(k)(r)}\bigr\rbrace_{k\in \mathbb{N}}$ is a basis for the
space $t_{p}^r$, and any $x\in t_{p}^r$ has a unique representation of the form
\begin{eqnarray}\label{eqb1}
x=\sum_{k}\lambda_{k}(r)b^{(k)}(r).
\end{eqnarray}

\begin{rem}

\begin{itemize}
  \item[(i)] It is well known that every Banach space $X$ with Schauder basis is seperable.
  \item[(ii)] Since the sequence space $\ell_{\infty}$ has no Schauder basis, Taylor seqeunce space $t_{\infty}^r$
has no Schuder basis. Also, the space $t_{\infty}^r$ is not seperable.
 \end{itemize}
\end{rem}

\begin{cor}
The space $t_{p}^{r}$ is seperable.
\end{cor}

In following theorem, we give some inclusion relations for the Taylor seqeunce spaces $t_{p}^r$ and $t_{\infty}^r$. These inclusion relations
can be proved as similar to  inclusion relations theorems in Altay et al.\cite{abm}, Mursaleen et al. \cite{mba} and Kirişci \cite{kiris2}.\\

\begin{thm}We have:
\begin{itemize}
 \item[(i)] The inclusions $\ell_{p}\subset t_{p}^{r}$ strictly holds for $1\leq p < \infty$.
 \item[(ii)] Neither of the spaces $t_{p}^{r}$ and $\ell_{\infty}$ includes the other one, where $1\leq p < \infty$.
 \item[(iii)] The space $t_{\infty}^{r}$ strictly includes both the space $\ell_{\infty}$ and the spaces $t_{p}^{r}$, where
$1\leq p < \infty$.
 \item[(iv)] If $0<s\leq r<1$, then $t_{p}^{r}\subset t_{p}^{s}$.
 \item[(v)] The inclusion $t_{p}^r\subset t_{s}^r$ holds whenever $1\leq p < s$.
\end{itemize}
\end{thm}

We assume througout that $p^{-1}+q^{-1}=1$ for $p,q\geq 1$ and denote the collection of all finite subsets on $\mathbb{N}$ by $\mathcal{F}$. Firstly, we give some results for use in proof of the Theorems \ref{alphathm}-\ref{gammathm}.\\

\begin{eqnarray}\label{90}
&&\sup_{K\in \mathcal{F}}\sum_{k}\left|\sum_{n\in K}a_{nk}\right|^{q}<\infty \quad (1<p\leq \infty),\\ \label{91}
&&\lim_{n\rightarrow \infty}a_{nk}=\alpha_{k} \quad (k\in \mathbb{N}),\\ \label{92}
&&\sup_{n\in \mathbb{N}}\sum_{k}|a_{nk}|^{q}<\infty \quad (1<p< \infty)
\end{eqnarray}

\begin{lem}\label{lemMTRX}
For the characterization of the class $(X:Y)$ with
$X=\{\ell_{p}\}$ and $Y=\{\ell_{1}, c, \ell_{\infty}\}$, we can give the necessary and sufficient
conditions from Table 1, where
\begin{center}
\begin{tabular}{|l | l | l |}
\hline \textbf{1.} (\ref{90}) & \textbf{2.} (\ref{91}), (\ref{92}) & \textbf{3.}  (\ref{92}) \\
\hline
\end{tabular}
\end{center}

\end{lem}

\begin{center}
\begin{tabular}{|c | c c c|}
\hline
To $\rightarrow$ & $\ell_{1}$ & $c$ & $\ell_{\infty}$\\ \hline
From $\downarrow$ &  &  &\\ \hline
$\ell_{p}$ & \textbf{1.} & \textbf{2.} & \textbf{3.}\\
\hline
\end{tabular}

\vspace{0.1cm}Table 1\\

\end{center}

Now, we may give the theorems determining the $\alpha-, \beta-$ and $\gamma-$duals of the Taylor sequence spaces
$t_{p}^{r}$ for $1\leq p\leq \infty$.\\

\begin{thm}\label{alphathm}
The $\alpha-$duals of the spaces $t_{p}^{r}=\alpha_{r}$ and $t_{1}^{r}=\alpha_{\infty}$, where
\begin{eqnarray*}
\alpha_{r}&=&\left\{a=(a_{k})\in \omega: \sup_{N,K\in \mathcal{F}}\sum_{n\in N}\left|\sum_{k\in K}c_{nk}^{r}\right|^{q}<\infty\right\}\\
\alpha_{\infty}&=&\left\{a=(a_{k})\in \omega: \sup_{n\in \mathbb{N}}\sum_{k}\left|c_{nk}^{r}\right|<\infty\right\}.
\end{eqnarray*}
where the matrix $C(r)=c_{nk}^{r}$ defined by
\begin{eqnarray}\label{mtrxalpha}
c_{nk}^r = \left\{ \begin{array}{ccl}
\binom{k}{n}(-r)^{k-n}(1-r)^{-(k+1)}a_{n}&, & (k\geq n),\\
0&, & (0\leq k < n)
\end{array} \right.
\end{eqnarray}

\end{thm}

\begin{proof}
We choose the sequence $a=(a_{k})\in \omega$. We can easily derive that with the (\ref{transseq}) that
\begin{eqnarray}\label{alphaeq1}
a_{n}x_{n}=\sum_{k=n}^{\infty}\binom{k}{n}(-r)^{k-n}(1-r)^{-(k+1)}a_{n}y_{k}=\left(C(r)y\right)_{n}, \quad (n\in \mathbb{N}).
\end{eqnarray}
for all $k,n\in \mathbb{N}$, where $C(r)=c_{nk}^{r}$ defined by (\ref{mtrxalpha}). It follows from (\ref{alphaeq1}) that $ax=(a_{n}x_{n})\in\ell_{1}$
whenever $x\in t_{p}^{r}$ or $t_{\infty}^{r}$ if and only if $Cy\in \ell_{1}$ whenever $y\in \ell_{p}$ or $\ell_{\infty}$. We obtain that $a\in (t_{p}^{r})^{\alpha}$ or
$a\in (t_{\infty}^{r})^{\alpha}$ whenever $x\in (t_{p}^{r})$ or $x\in (t_{\infty}^{r})$ if and only if $C\in (\ell_{p}:\ell_{1})$ or $C\in (\ell_{\infty}:\ell_{1})$.
Therefore, we get by Lemma \ref{lemMTRX} with the matrix $C$ instead of $A$ that $a\in (t_{p}^{r})^{\alpha}$ or
$a\in (t_{\infty}^{r})^{\alpha}$ if and only if $\sup_{n\in K}\sum_{k}|c_{nk}^{r}|<\infty$. This gives us the result that
$(t_{p}^{r})^{\alpha}=\alpha_{r}$ and $(t_{\infty}^{r})^{\alpha}=\alpha_{r}$.
\end{proof}

\begin{thm}\label{betathm}
The matrix $D(r)=(d_{nk}^{r})$ is defined by
\begin{eqnarray}\label{betamtrx}
d_{nk}^r = \left\{ \begin{array}{ccl}
\sum_{k=0}^{n}\binom{n}{k}(-r)^{n-k}(1-r)^{-(n+1)}a_{k}&, & (0\leq k \leq n),\\
0&, & (k>n)
\end{array} \right.
\end{eqnarray}
for all $k,n\in\mathbb{N}$. Then, $(t_{1}^{r})^{\beta}=\beta_{2}\cap \beta_{4}$,
$(t_{p}^{r})^{\beta}=\beta_{1}\cap \beta_{2}$ and $(t_{\infty}^{r})^{\beta}=\beta_{2}\cap \beta_{3}$, where
\begin{eqnarray*}
\beta_{1}&=&\left\{a=(a_{k})\in \omega: \sup_{ n\in \mathbb{N}}\sum_{k}|d_{nk}^{r}|^{q}<\infty\right\} \quad (1<p<\infty),\\
\beta_{2}&=&\left\{a=(a_{k})\in \omega: \lim_{n\rightarrow\infty}d_{nk}^{r}\quad \textrm{exists for each}\quad k\in \mathbb{N}\right\},\\
\beta_{3}&=&\left\{a=(a_{k})\in \omega: \lim_{n\rightarrow\infty}\sum_{k}|d_{nk}^{r}| \quad \textrm{exists}\right\},\\
\beta_{4}&=&\left\{a=(a_{k})\in \omega: \sup_{ k,n\in \mathbb{N}}|d_{nk}^{r}|<\infty\right\}.
\end{eqnarray*}
\end{thm}

\begin{proof}
Consider the equation
\begin{eqnarray}\label{betaeq}
\sum_{k=0}^{n}a_{k}x_{k}=\sum_{k=0}^{n}\left[\sum_{k=j}^{\infty} \binom{k}{j}(-r)^{k-j}(1-r)^{-(k+1)}y_{k}\right]a_{k}
\end{eqnarray}
\begin{eqnarray*}
\quad \quad \quad \quad \quad \quad \quad \quad
=\sum_{k=0}^{n}\left[\sum_{j=0}^{k} \binom{k}{j}(-r)^{k-j}(1-r)^{-(k+1)}a_{j}\right]y_{k}=(D(r)y)_{n}
\end{eqnarray*}
Then, we deduce from Lemma \ref{lemMTRX} with (\ref{betaeq}) that $ax=(a_{k}x_{k})\in cs$ whenever
$x\in t_{p}^{r}$ if and only if $D(r)y\in c$ whenever $y\in \ell_{p}$. That is to say that
$a=(a_{k})\in (t_{p}^{r})^{p}$ if and only if $D(r)\in (\ell_{p}:c)$. Therefore, we derive from (\ref{91}) and (\ref{92}),
we conclude that $\lim_{n\rightarrow\infty}d_{nk}^{r}$ exists and $\sup_{ n\in \mathbb{N}}\sum_{k}|d_{nk}^{r}|<\infty$
which shows that $(t_{p}^{r})^{\beta}=\beta_{1}\cap \beta_{2}$.\\

In a similar way, we can prove that $(t_{1}^{r})^{\beta}=\beta_{2}\cap \beta_{4}$, and $(t_{\infty}^{r})^{\beta}=\beta_{2}\cap \beta_{3}$.
\end{proof}

\begin{thm}\label{gammathm}
$(t_{1}^{r})^{\gamma}=\beta_{4}$, $(t_{p}^{r})^{\gamma}=\beta_{1} \quad (1<p<\infty)$ and $(t_{\infty}^{r})^{\gamma}=\beta_{5}$,
where
\begin{eqnarray*}
\beta_{5}=\left\{a=(a_{k})\in \omega: \sup_{ n\in \mathbb{N}}\sum_{k}|d_{nk}^{r}|<\infty\right\} \quad (1<p<\infty)
\end{eqnarray*}

\end{thm}
\begin{proof}
Let $a=(a_{k})\in\beta_{1}$ and $x=(x_{k})\in t_{p}^{r}$. Then, we obtain by applying the H\"{o}lder's inequality that
\begin{eqnarray*}
\left|\sum_{k=0}^{n}a_{k}x_{k}\right|&=&\left|\sum_{k=0}^{n}\left[\sum_{k=j}^{\infty} \binom{k}{j}(-r)^{k-j}(1-r)^{-(k+1)}y_{k}\right]a_{k}\right|\\
&=&\left|\sum_{k=0}^{n}d_{nk}^{r}y_{k}\right|\leq \left(\sum_{k=0}^{n}|d_{nk}^{r}|^{q}\right)^{1/q}\left(\sum_{k=0}^{n}|y_{k}|^{p}\right)^{1/p}
\end{eqnarray*}
where $d_{nk}^{r}$ is defined by (\ref{betamtrx}). Taking supremum over $n\in \mathbb{N}$, we have
\begin{eqnarray*}
\sup_{n\in \mathbb{N}}\left|\sum_{k=0}^{n}a_{k}x_{k}\right|
&\leq& \sup_{n\in \mathbb{N}}\left[\left(\sum_{k=0}^{n}|d_{nk}^{r}|^{q}\right)^{1/q}\left(\sum_{k=0}^{n}|y_{k}|^{p}\right)^{1/p}\right]\\
&\leq& \|y\|_{\ell_{p}}\left(\sup_{n\in \mathbb{N}}\sum_{k=0}^{n}|d_{nk}^{r}|^{q}\right)^{1/q}\leq \infty.
\end{eqnarray*}
This means that $a=(a_{k})\in (t_{p}^{r})^{\gamma}$. Hence $\beta_{1}\subset (t_{p}^{r})^{\gamma}$.\\

Conversely, let $a=(a_{k})\in (t_{p}^{r})^{\gamma}$ and $x=(x_{k})\in t_{p}^{r}$. Then one can easily see that
$\{\sum_{k=0}^{n}d_{nk}^{r}y_{k}\}\in \ell_{\infty}$, for all $n\in\mathbb{N}$, whenever $(a_{k}x_{k})\in bs$. This implies that the triangle $D(r)=(d_{nk}^{r})$
defined by (\ref{betamtrx}), is in the class $(\ell_{p}:\ell_{\infty})$. Hence, the condition $\sup_{ n\in \mathbb{N}}\sum_{k}|d_{nk}^{r}|<\infty$
is satisfied, which implies that $a=(a_{k})\in \beta_{1}$. That is, $(t_{p}^{r})^{\gamma} \subset \beta_{1}$. Therefore, we establish that the
$\gamma-$dual of $t_{p}^{r}$ is the set $\beta_{1}$.\\

With the same idea, it can be proved that $(t_{1}^{r})^{\gamma}=\beta_{4}$ and $(t_{\infty}^{r})^{\gamma}=\beta_{5}$.
\end{proof}

\begin{lem}\cite[Lemma 5.3]{AB2}\label{mtrxtool0}
Let $X, Y$ be any two sequence spaces, $A$ be an infinite matrix and $U$ a triangle matrix matrix.Then, $A\in (X: Y_{U})$ if and only if $UA\in (X:Y)$.
\end{lem}

\begin{lem}\cite[Theorem 3.1]{AB}\label{mtrxtool}
$B^{U}=(b_{nk})$ be defined via a sequence $a=(a_{k})\in\omega$ and inverse of the triangle matrix $U=(u_{nk})$ by
\begin{eqnarray*}
b_{nk}=\sum_{j=k}^na_{j}v_{jk}
\end{eqnarray*}
for all $k,n\in\mathbb{N}$. Then,
\begin{eqnarray*}
\lambda_{U}^{\beta}=\{a=(a_{k})\in\omega: B^{U}\in(\lambda:c)\}
\end{eqnarray*}
and
\begin{eqnarray*}
\lambda_{U}^{\gamma}=\{a=(a_{k})\in\omega: B^{U}\in(\lambda:\ell_{\infty})\}.
\end{eqnarray*}
\end{lem}

In what follows, for brevity, we write,
\begin{eqnarray*}
\overline{a}_{nk}:=\sum_{k=0}^{n}\binom{n}{k}(-r)^{n-k}(1-r)^{-(n+1)}a_{nk}
\end{eqnarray*}
for all $k,m,n\in\mathbb{N}$. From Lemma \ref{mtrxtool0} and Lemma \ref{mtrxtool}, we have:

\begin{thm}\label{th41}
Suppose that the entries of the infinite matrices $A=(a_{nk})$ and $E=(e_{nk})$ are connected with the relation
\begin{eqnarray}\label{mtrxeq1}
e_{nk}=\overline{a}_{nk}
\end{eqnarray}
for all $k,n\in \mathbb{N}$ and $\mu$ be any given sequence space. Then,
\begin{itemize}
  \item[(i)] $A \in (t_{p}^{r}:\mu)$ if and only if $\{a_{nk}\}_{k\in\mathbb{N}} \in \{t_{p}^{r}\}^{\beta}$
  for all $n\in \mathbb{N}$ and $E\in (\ell_{p}:\mu)$.
  \item[(ii)] $A \in (t_{\infty}^{r}:\mu)$ if and only if $\{a_{nk}\}_{k\in\mathbb{N}} \in \{t_{\infty}^{r}\}^{\beta}$
  for all $n\in \mathbb{N}$ and $E\in (\ell_{\infty}:\mu)$.
\end{itemize}
\end{thm}

\begin{thm}\label{th42}
 Suppose that the elements of the infinite matrices $A=(a_{nk})$
and $B=(b_{nk})$ are connected with the relation
\begin{eqnarray}\label{bnk}
b_{nk}:=\sum_{j=n}^{\infty}\binom{j}{n}(1-r)^{n+1}r^{j-n}a_{jk}~\text{ for all }~k,n\in \mathbb{N}.
\end{eqnarray}
Let $\mu $ be any given sequence space. Then,
\begin{itemize}
  \item[(i)] $A \in (\mu : t_{p}^{r})$ if and only if $B\in (\mu: \ell_{p})$.
  \item[(ii)] $A \in (\mu : t_{\infty}^{r})$ if and only if $B\in (\mu: \ell_{\infty})$.
\end{itemize}

\end{thm}

The folowing results were taken from Stieglitz and Tietz \cite{ST}:
\begin{eqnarray}\label{3}
&&\sum_{n}a_{nk} \quad \textrm{convergent for all $k$},\\ \label{5}
&&\sup_{k,m}\left|\sum_{n=0}^{m}a_{nk}\right|<\infty,\\ \label{6}
&&\sup_{n,k}|a_{nk}|<\infty,\\ \label{7}
&&\sum_{n}a_{nk}=0  \quad \textrm{for all $k$},\\ \label{10}
&&\lim_{m}\sum_{k} \left|\sum_{n=0}^{m}a_{nk}\right|=\sum_{k}\left|\sum_{n}a_{nk}\right|,\\ \label{11}
&&\lim_{m}\sum_{k} \left|\sum_{n=0}^{m}a_{nk}\right|=0\\ \label{14}
&&\sup_{m}\sum_{k} \left|\sum_{n=0}^{m}a_{nk}\right|^{q}<\infty\\ \label{15}
&&\lim_{k}a_{nk}=0 \quad \textrm{for all $n$}\\ \label{16}
&&\sup_{n}\sum_{k}\left|a_{nk}-a_{n,k+1}\right|<\infty\\ \label{17}
&&\sup_{K\in\mathcal{F}}\sum_{n} \left|\sum_{k\in K}(a_{nk}-a_{n,k+1})\right|^{p}<\infty\\ \label{18}
&&\sup_{K,N\in\mathcal{F}}\left|\sum_{n\in N} \sum_{k\in K}(a_{nk}-a_{n,k+1})\right|<\infty\\ \label{19}
&&\sup_{n}\left|\lim_{k}a_{nk}\right|<\infty\\ \label{20}
&&\sup_{K\in\mathcal{F}}\sum_{n} \left|\sum_{k\in K}(a_{nk}-a_{n,k-1})\right|^{p}<\infty\\ \label{21}
&&\sup_{K,N\in\mathcal{F}}\left|\sum_{n\in N} \sum_{k\in K}(a_{nk}-a_{n,k-1})\right|<\infty\\ \label{5A}
&&\lim_{n}\sum_{k}|a_{nk}|=\sum_{k}\left|\lim_{n}a_{nk}\right|\\ \label{5B}
&&\lim_{n}\sum_{k}|a_{nk}|=0\\ \label{5C}
&&\sup_{K\in\mathcal{F}}\sum_{n} \left|\sum_{k\in K}a_{nk}\right|^{p}<\infty
\end{eqnarray}

\begin{lem}\label{lemMTRX2}
For the characterization of the class $(X:Y)$ for
$\{\ell_{\infty}, \ell_{p}, \ell_{1}, c, c_{0}, bs, cs, c_{0}s\}$, we can give the necessary and sufficient
conditions from Table 2, Table 3, Table 4 and Table 5 where\\
\begin{center}
\begin{tabular}{|l | l | l | l |}
\hline \textbf{4.} (\ref{92}) for $q=1$ & \textbf{5.} (\ref{91}), (\ref{5A}) & \textbf{6.}  (\ref{5B}) & \textbf{7.} (\ref{91}) for $\alpha_{k}=0, $ (\ref{92})  \\
\hline \textbf{8.} (\ref{6}) &  \textbf{9.} (\ref{91}), (\ref{6}) & \textbf{10.} (\ref{91}) for $\alpha_{k}=0, $ (\ref{6}) & \textbf{11.} (\ref{5C})\\
\hline  \textbf{12.} (\ref{5C}) for $p=1$ & \textbf{13.} (\ref{14}) & \textbf{14.} (\ref{3}), (\ref{14})  & \textbf{15.} (\ref{7}), (\ref{14}) \\
\hline \textbf{16.} (\ref{5})  & \textbf{17.} (\ref{3}), (\ref{5})& \textbf{18.}  (\ref{5}), (\ref{7})  & \textbf{19.}  (\ref{14}) for $q=1$ \\
\hline  \textbf{20.} (\ref{3})   & \textbf{21.} (\ref{11}) & \textbf{22.}  (\ref{15}), (\ref{16}) & \textbf{23.} (\ref{15}), (\ref{17})  \\
\hline  \textbf{24.} (\ref{15}), (\ref{18})  &  \textbf{25.}   (\ref{16}), (\ref{19}) & \textbf{26.} (\ref{20})   & \textbf{27.} (\ref{21}) \\
\hline \textbf{28.} (\ref{16})   & \textbf{29.} (\ref{17})  & \textbf{30.} (\ref{18}) & \\
\hline
\end{tabular}
\end{center}
\end{lem}

\begin{center}
\begin{tabular}{|c | c c c|}
\hline
To $\rightarrow$ & $\ell_{\infty}$ & $c$ & $c_{0}$\\ \hline
From $\downarrow$ &  &  &\\ \hline
$\ell_{\infty}$ & \textbf{4.} & \textbf{5.} & \textbf{6.}\\ \hline
$\ell_{p}$ & \textbf{3.} & \textbf{2.} & \textbf{7.} \\ \hline
$\ell_{1}$ & \textbf{8.} & \textbf{9.} & \textbf{10.}\\
\hline
\end{tabular}

\vspace{0.1cm}Table 2\\

\end{center}

\begin{center}
\begin{tabular}{|c | c c c|}
\hline
To $\rightarrow$ & $\ell_{\infty}$ & $\ell_{p}$ & $\ell_{1}$\\ \hline
From $\downarrow$ &  &  &\\ \hline
$\ell_{\infty}$ & \textbf{4.} & \textbf{11.} & \textbf{12.}\\ \hline
$c$ & \textbf{4.} & \textbf{11.} & \textbf{12.} \\ \hline
$c_{0}$ & \textbf{4.} & \textbf{11.} & \textbf{12.}\\
\hline
\end{tabular}

\vspace{0.1cm}Table 3\\

\end{center}

\begin{center}
\begin{tabular}{|c | c c c|}
\hline
To $\rightarrow$ & $bs$ & $cs$ & $c_{0}s$\\ \hline
From $\downarrow$ &  &  &\\ \hline
$\ell_{p}$ & \textbf{13.} & \textbf{14.} & \textbf{15.}\\ \hline
$\ell_{1}$ & \textbf{16.} & \textbf{17.} & \textbf{18.} \\ \hline
$\ell_{\infty}$ & \textbf{19.} & \textbf{20.} & \textbf{21.}\\
\hline
\end{tabular}

\vspace{0.1cm}Table 4\\

\end{center}

\begin{center}
\begin{tabular}{|c | c c c|}
\hline
To $\rightarrow$ & $\ell_{\infty}$ & $\ell_{p}$ & $\ell_{1}$\\ \hline
From $\downarrow$ &  &  &\\ \hline
$bs$ & \textbf{22.} & \textbf{23.} & \textbf{24.}\\ \hline
$cs$ & \textbf{25.} & \textbf{26.} & \textbf{27.} \\ \hline
$c_{0}s$ & \textbf{28.} & \textbf{29.} & \textbf{30.}\\
\hline
\end{tabular}

\vspace{0.1cm}Table 5\\

\end{center}

Now, using the Theorem \ref{th41}, Theorem \ref{th42} and Lemma \ref{lemMTRX2}, we can give the some results:

Let $A=(a_{nk})$ be an infinite matrix. Consider the $X\in \{\ell_{\infty}, c, c_{0}, bs, cs, c_{0}s\}$.\\

$A\in (t_{p}^{r}:X)$ if and only if $\{a_{nk}\}_{k\in \mathbb{N}}\in \{t_{p}^{r}\}^\beta$ for all
$n\in\mathbb{N}$ and the conditions \textbf{2.}, \textbf{3.}, \textbf{7.}, \textbf{13.}, \textbf{14.}, \textbf{15.} in Table 2 and Table 4
hold with $\overline{a}_{nk}$ instead of $a_{nk}$.\\

$A\in (t_{1}^{r}:X)$ if and only if $\{a_{nk}\}_{k\in \mathbb{N}}\in \{t_{1}^{r}\}^\beta$ for all
$n\in\mathbb{N}$ and the conditions \textbf{8.}, \textbf{9.}, \textbf{10.}, \textbf{16.}, \textbf{17.}, \textbf{18.} in Table 2 and Table 4
hold with $\overline{a}_{nk}$ instead of $a_{nk}$.\\

$A\in (t_{\infty}^{r}:X)$ if and only if $\{a_{nk}\}_{k\in \mathbb{N}}\in \{t_{\infty}^{r}\}^\beta$ for all
$n\in\mathbb{N}$ and the conditions \textbf{4.}, \textbf{5.}, \textbf{6.}, \textbf{19.}, \textbf{20.}, \textbf{21.} in Table 2 and Table 4
hold with $\overline{a}_{nk}$ instead of $a_{nk}$.\\

The relation $b_{nk}$ be defined by (\ref{bnk}).\\

$A=(a_{nk})\in (X:t_{\infty}^{r})$ if and only if the conditions \textbf{4.}, \textbf{22.}, \textbf{25.}, \textbf{28.} in Table 3 and Table 5
hold with $b_{nk}$ instead of $a_{nk}$.\\

$A=(a_{nk})\in (X:t_{p}^{r})$ if and only if the conditions \textbf{11.}, \textbf{23.}, \textbf{26.}, \textbf{29.} in Table 3 and Table 5
hold with $b_{nk}$ instead of $a_{nk}$.\\

$A=(a_{nk})\in (X:t_{1}^{r})$ if and only if the conditions \textbf{12.}, \textbf{24.}, \textbf{27.}, \textbf{30.} in Table 3 and Table 5
hold with $b_{nk}$ instead of $a_{nk}$.\\

Thus, we have characterized the matrix classes between the new space and the classical spaces. In above results, some examples are given with the
different spaces wihch are obtained from the triangular matrices(cf. \cite{B4}) instead of $X$.

\section{Some Geometric Properties of the Space $t^{r}_{p}$}

In this section, we examine some geometric properties of the space $t^{r}_{p}$.
First, we define some geometric properties of the spaces. Let $(X,\|.\|)$ be a normed space and let $S(x)$ and $B(x)$ be the unit sphere and unit ball of $X$, respectively. Consider \textit{Clarkson's modulus of convexity} (see \cite{clark1,clark2}) defined by
\begin{eqnarray*}
\delta_{X}(\theta)=\inf\left\{1-\frac{\| x-y\|}{2}; ~~ x,y\in S(x),~\| x-y\|=\theta \right\},
\end{eqnarray*}
where $0\leq\theta\leq 2$. The inequality $\delta_{X}(\theta)>0$ for all $\theta\in [0,2]$ characterizes the uniformly convex
spaces.
 In \cite{gur}, \textit{Gurarii's modulus of convexity} is defined by
\begin{eqnarray*}
\beta_{X}(\theta)=\inf\left\{1-\inf_{\alpha\in [0,1]}\| \alpha x+(1-\alpha)y\|; ~~ x,y\in S(x),~\| x-y\|=\theta\right\},
\end{eqnarray*}
where $0\leq\theta\leq 2$. It is easily shown that $\delta_{X}(\theta)\leq \beta_{X}(\theta) \leq 2\delta_{X}(\theta)$ for any $0\leq\theta\leq 2$. Further, if  $0<\beta_{X}(\theta)<1$, then $X$ is uniformly convex, and if $\beta_{X}(\theta)<1$, then $X$ is strictly convex.

A Banach space $X$ is said to have the Banach-Saks property if every bounded sequence $(x_{n})$ in $X$ admits a sequence $(z_{n})$ such that the sequence
$t_{k}(z)$ is convergent in norm in  $X$ (see \cite{diz}), where
\begin{eqnarray*}
t_{k}(z)=\frac{1}{k+1}(z_{0}+z_{1}+\cdots+z_{k}) ~\textrm{ for all }~k\in\mathbb{N}.
\end{eqnarray*}
Let $1<p<\infty$. A Banach space is said to have the Banach-Saks type $p$ if every weakly null sequence has a subsequence $(x_{k})$ such that for some $C>0$,
\begin{eqnarray*}
\left\|\sum_{k=0}^{n} x_{k}\right \|< C(n+1)^{1/p}.
\end{eqnarray*}

A Banach space $X$ is said to have the \textit{weak Banach-Saks property} whenever given any weakly null sequence $(x_{n})$ in $X$ and there exists a subsequence $(z_{n})$ of $(x_{n})$ such that the sequence $\{t_{k}(z)\}$ strongly converges to zero.

In \cite{gar}, Garc\'{\i}a-Falset introduced the following coefficient:
\begin{eqnarray*}
R(X)=\sup\left\{\liminf_{n\to\infty} \| x_{n}+x\|; ~~ (x_{n})\subset B(x),~x_{n}\stackrel{w}{\to} 0\right\}.
\end{eqnarray*}

\begin{rem}\label{rem1}\cite{gar1}
A Banach space $X$ with $R(X)<2$ has a \textit{weak fixed point property}.
\end{rem}
\begin{thm}\label{lem2.14}
The space $t^{r}_{p}$ has Banach-Saks type $p$.
\end{thm}
\begin{proof}
Let $(\varepsilon_{n})$ be a sequence of positive numbers for which $\sum_{n=1}^{\infty}\varepsilon_{n}\leq 1/2$. Let $(x_{n})$ be a weakly null sequence in
$B(t^{p}(r))$. Set $s_{0}=x_{0}$ and $s_{1}=x_{n_{1}}=x_{1}$. Then, there exists $t_{1}\in\mathbb{N}$ such that
\begin{eqnarray*}
\left\|\sum_{i=t_{1}+1}^{\infty} s_{1}(i)e^{(i)}\right \|_{t^{r}_{p}}<\varepsilon_{1}.
\end{eqnarray*}
The assumption "$(x_{n})$ be a weakly null sequence" implies that $x_{n}\to 0$ with respect to the coordinatewise, there exists an  $n_{2}\in\mathbb{N}$ such that
\begin{eqnarray*}
\left\|\sum_{i=0}^{t_{1}} x_{n}(i)e^{(i)}\right \|_{t^{r}_{p}}<\varepsilon_{1},
\end{eqnarray*}
where $n\geq n_{2}$. Set $s_{2}=x_{n_{2}}$. Then, there exists $t_{2}>t_{1}$ such that

\begin{eqnarray*}
\left\|\sum_{i=t_{2}+1}^{\infty} s_{2}(i)e^{(i)}\right \|_{t^{r}_{p}}<\varepsilon_{2}.
\end{eqnarray*}
By using the fact that $x_{n}\to 0$ with respect to coordinatewise, there exists $n_{3}>n_{2}$  such that
\begin{eqnarray*}
\left\|\sum_{i=0}^{t_{2}} x_{n}(i)e^{(i)}\right \|_{t^{r}_{p}}<\varepsilon_{2},
\end{eqnarray*}
where $n\geq n_{3}$. If we continue this process, we can find two increasing sequences $(t_{i})$ and $(n_{i})$ of natural numbers such that
\begin{eqnarray*}
\left\|\sum_{i=0}^{t_{j}} x_{n}(i)e^{(i)}\right \|_{t^{r}_{p}}<\varepsilon_{j}
\end{eqnarray*}
for each $n\geq n_{j+1}$ and
\begin{eqnarray*}
\left\|\sum_{i=t_{j}+1}^{\infty} s_{j}(i)e^{(i)}\right \|_{t^{r}_{p}}<\varepsilon_{j},
\end{eqnarray*}
where $s_{j}=x_{n_{j}}$. Hence,
\begin{eqnarray*}
\left\|\sum_{j=0}^{n} s_{j}\right \|_{t^{r}_{p}}&=&\left\|\sum_{j=0}^{n}\left(\sum_{i=0}^{t_{j}-1}s_{j}(i)e^{(i)}+ \sum_{i=t_{j}-1}^{t_{j}}s_{j}(i)e^{(i)}+ \sum_{i=t_{j}+1}^{\infty}s_{j}(i)e^{(i)}\right)\right \|_{t^{r}_{p}}\\
&\leq & \left\|\sum_{j=0}^{n}\sum_{i=0}^{t_{j}-1}s_{j}(i)e^{(i)}\right \|_{t^{p}(r)}+ \left\|\sum_{j=0}^{n}\sum_{i=t_{j}-1}^{t_{j}}s_{j}(i)e^{(i)}\right \|_{t^{r}_{p}} \\
&+&\left\|\sum_{j=0}^{n}\sum_{i=t_{j}+1}^{\infty}s_{j}(i)e^{(i)}\right \|_{t^{r}_{p}}\\
&\leq& \left\|\sum_{j=0}^{n}\left( \sum_{i=t_{j-1}+1}^{t_{j}}s_{j}(i)e^{(i)}\right)\right \|_{t^{r}_{p}} + 2 \sum_{j=0}^{n}\varepsilon_{j}.
\end{eqnarray*}
In other respects, one can see that $\| x \|_{t^{r}_{p}}<1$. Thus, $\|x\|^{p}_{t^{r}_{p}}<1$ and we get
\begin{eqnarray*}
\left\|\sum_{j=0}^{n}\sum_{i=t_{j-1}+1}^{t_{j}}s_{j}(i)e^{(i)}\right\|^{p}_{t^{r}_{p}}&=& \sum_{j=0}^{n}\sum_{i=t_{j}-1}^{t_{j}}\left|\sum_{k=i}^{\infty}{k\choose i}(1-r)^{i+1}r^{k-i}x_{j}(k)\right|^{p}\\
&\leq&\sum_{j=0}^{n}\sum_{i=0}^{\infty}\left|\sum_{k=i}^{\infty}{k\choose i}(1-r)^{i+1}r^{k-i}x_{j}(k)\right|^{p}\\
&\leq& (n+1).
\end{eqnarray*}
Hence, we get
\begin{eqnarray*}
\left\|\sum_{j=0}^{n}\sum_{i=t_{j-1}+1}^{t_{j}}s_{j}(i)e^{(i)}\right\|_{t^{r}_{p}}\leq (n+1)^{1/p}.
\end{eqnarray*}
By using the inequality $1\leq (n+1)^{1/p}$ for all $n\in\mathbb{N}$ and $1\leq p <\infty$, we obtain
\begin{eqnarray*}
\left\|\sum_{j=0}^{n} s_{j}\right \|_{t^{r}_{p}}\leq (n+1)^{1/p}+1 \leq 2(n+1)^{1/p}.
\end{eqnarray*}
Thus, the space $t^{r}_{p}$ has Banach-Saks type $p$.
\end{proof}

\begin{rem}\label{rem2}
Note that $R(t^{r}_{p})=R(\ell_{p})=2^{1/p}$, since $t^{r}_{p}$ is an isomorphic copy of the $\ell_{p}$.
\end{rem}
Hence, by Remarks \ref{rem1} and \ref{rem2}, we have the following:
\begin{cor}
Let $1<p<\infty$. Then, the space $t^{r}_{p}$ has the weak fixed point property.
\end{cor}
\begin{thm}
Gurarii's modulus of convexity for the normed space $t^{r}_{p}$ is
\begin{eqnarray*}
\beta_{t^{r}_{p}}(\theta)\leq 1-\left[1-\left(\frac{\theta}{2}\right)^{p}  \right]^{1/p},
\end{eqnarray*}
where $0\leq \theta\leq 2$.
\end{thm}
\begin{proof}
Let $x\in t^{r}_{p}$. Then, we obtain
\begin{eqnarray*}
\| x\|_{t^{r}_{p}}=\| T(r)x\|_{p}=\left[\sum_{n}\left|(T(r)x)_{n}\right|^{p}\right]^{1/p}.
\end{eqnarray*}
Let $0\leq \theta\leq 2$ and take into consideration this sequences
\begin{eqnarray*}
x=(x_{k})&=&\left\{T[-r/(1-r)]\left[\left[1-\left(\frac{\theta}{2}\right)^{p} \right]^{1/p}\right],T[-r/(1-r)]\left( \frac{\theta}{2} \right),0,0,\ldots \right\},\\
y=(y_{k})&=&\left\{T[-r/(1-r)]\left[\left[1-\left(\frac{\theta}{2}\right)^{p} \right]^{1/p}\right],T[-r/(1-r)]\left(-\frac{\theta}{2}\right),0,0,\ldots\right\}.
\end{eqnarray*}
Because of $a_{k}=[T(r)x]_{k}$ and $b_{k}=[T(r)y]_{k}$, we get
\begin{eqnarray*}
a=(a_{k})&=&\left\{\left[1-\left(\frac{\theta}{2}\right)^{p}\right]^{1/p},\frac{\theta}{2},0,0,\ldots\right\},\\
b=(b_{k})&=&\left\{\left[1-\left(\frac{\theta}{2}\right)^{p}\right]^{1/p} , -\frac{\theta}{2},0,0,\ldots\right\}.
\end{eqnarray*}
By using the sequences $x=(x_{k})$ and $y=(y_{k})$, we obtain following equalities
\begin{eqnarray*}
\| x \|^{p}_{t^{r}_{p}}= \|T(r)x\|_{p}&=& \left|\left[1-\left(\frac{\theta}{2}\right)^{p} \right]^{1/p}\right|^{p}+\left|-\frac{\theta}{2}\right|^{p} \\
&=&1-\left(\frac{\theta}{2}\right)^{p}+\left(\frac{\theta}{2}\right)^{p}=1,\\
\end{eqnarray*}
\begin{eqnarray*}
\| y \|^{p}_{t^{r}_{p}}=\|T(r)y\|_{p}&=& \left|\left[1-\left(\frac{\theta}{2}\right)^{p}\right]^{1/p}\right|^{p}+\left|-\frac{\theta}{2}\right|^{p} \\
&=&1-\left(\frac{\theta}{2}\right)^{p}+\left(\frac{\theta}{2}\right)^{p}=1\\
\end{eqnarray*}
and
\begin{eqnarray*}
\|x-y\|^{p}_{t^{r}_{p}}=\|T(r)x-T(r)y\|_{p}&=& \left\{\left|\left[1-\left(\frac{\theta}{2}\right)^{p}\right]^{1/p}-\left[1-\left(\frac{\theta}{2}\right)^{p}\right]^{1/p}\right|^{p}+\left|\frac{\theta}{2}-\left(-\frac{\theta}{2}\right)\right|^{p}\right\}^{1/p} =\theta.
\end{eqnarray*}
For $0\leq\alpha\leq 1$
\begin{eqnarray*}
\inf_{\alpha\in[0,1]}\| \alpha x+(1-\alpha)y\|_{t^{r}_{p}}&=& \inf_{\alpha\in[0,1]}\| \alpha T(r)x+(1-\alpha)T(r)y\|_{p}\\
&=& \inf_{\alpha\in[0,1]}\left\{\left|\alpha\left[1-\left(\frac{\theta}{2}\right)^{p} \right]^{1/p}+(1-\alpha)\left[1-\left(\frac{\theta}{2}\right)^{p} \right]^{1/p}\right|^{p}+\left|\alpha\left(\frac{\theta}{2}\right)^{p}+(1-\alpha)\left(-\frac{\theta}{2}\right)\right|^{p}\right\}^{1/p}\\
&=&\inf_{\alpha\in[0,1]}\left[ 1-\left(\frac{\theta}{2}\right)^{p}+(2\alpha-1)^{p}\left(\frac{\theta}{2}\right)^{p} \right]^{1/p}\\
&=&\left[1-\left(\frac{\theta}{2}\right)^{p}\right]^{1/p}.
\end{eqnarray*}
Therefore, for $1\leq p<\infty$, we have
\begin{eqnarray*}
\beta_{t^{r}_{p}}(\theta)\leq 1-\left[1-\left(\frac{\theta}{2}\right)^{p}  \right]^{1/p}.
\end{eqnarray*}
The completes the proof.
\end{proof}
\begin{cor} The following statements hold:
\begin{enumerate}
\item[(i)] For $\theta >2$, $\beta_{t^{r}_{p}}(\theta)=1$. Thus, $t^{r}_{p}$ is strictly convex.
\item[(ii)] For $0<\theta \leq 2$, $\beta_{t^{r}_{p}}(\theta)\leq 1$. Thus, $t^{r}_{p}$ is uniformly convex.
\end{enumerate}
\end{cor}
\begin{cor}
 For $\alpha =1/2$, $\beta_{t^{r}_{p}}(\theta)=\delta_{t^{r}_{p}}(\theta)$.
\end{cor}

\section{Conclusion}
The purpose of this paper is twohold. In section 2, it was obtained the new matrix domain with the well known classical sequence spaces and an infinite matrix.
We use the Taylor method as an infinite matrix for the matrix domain, in this study. The Taylor method which known then as the circle method of order $r$ (0 < r < 1).\\

It is well known that every Banach space is isomorphic to a subspace of a Banach space and every separable reflexive Banach space
is isomorphic to a subspace of a separable reflexive space with Schauder basis. In section 2, we construct new sequence spaces,
which naturally emerge from the classical sequence spaces and Taylor transform. The new space $t^{r}_{p}$ obtained as domain of
Taylor upper triangle matrix in the space $\ell_{p}$ $(1\leq p \leq \infty)$. In this paper, since the new space $t^{r}_{p}$ is isomorphic copy to the space $\ell_{p}$, the space $t^{r}_{p}$ has a Schauder basis and also some topological properties of new constructed space are similar to the classical space $\ell_{p}$. Then, some basic results of the space $t^{r}_{p}$ such as absolute property, $BK-$space, isomorphism, Schauder basis, separability and some inclusion relations are given but not proved because of
natural continuation of the space $\ell_{p}$. Therefore, we show that the space $t^{r}_{p}$ is a Hilbert space for $p=2$. Also, the proofs of dual spaces theorems are given which are important. The characterized matrix classes between the new space and the classical spaces are shown in the table form. Further,
based on the tables, some results of matrix classes are given as examples.\\

In section 3, we focus on some important results which are geometric properties of Banach spaces. In Banach Space Theory, geometric properties
are play a crucial role. As is well known, among all infinite dimensional Banach spaces, Hilbert spaces have the nicest geometric properties.
We present the some geometric properties of the space  $t^{r}_{p}$ such as Clarkson's modulus of convexity, Gurarii's modulus of convexity, Banach-Saks
property, weak Banach-Saks property, weak fixed point property, strictly convexity, uniform convexity.\\

The geometric properties in this study will form the basis for future works. Then, the aim of this
paper is to present an in-depth and  up to date coverage of the main ideas, concepts and most important results related to
sequence spaces and geometric properties of Banach Space.

\section*{Acknowledgements}
 This work was supported by Scientific Projects Coordination Unit of Istanbul University. Project number 20969.

\section*{Competing interests}
The author declares that there is no conflict of interests regarding the publication of this article.

\end{document}